# On Hamilton Decompositions

Dhananjay P. Mehendale
Sir Parashurambhau College, Tilak Road, Pune-411030,
India

## Abstract

P. J. Kelly conjectured in 1968 that every diregular tournament on ($2n+1$) points can be decomposed in directed Hamilton circuits [1]. We define so called leading diregular tournament on ($2n+1$) points and show that it can be decomposed in directed Hamilton circuits when ($2n+1$) is a prime number. When ($2n+1$) is not a prime number this method does not work and we will need to devise some another method. We also propose a general method to find Hamilton decomposition of certain tournament for all sizes.

**1. Preliminaries:** We begin with few definitions: Let $G$ be a digraph with vertex set, $V(G) = \{v_1, v_2, \cdots, v_n\}$, and the set of directed edges, $E(G)$, and let $\delta^-(v_i)$ and $\delta^+(v_i)$ be the in-degree and out-degree of vertex $v_i$ respectively. We associate an **ordered pair** ($m_i, n_i$) with each vertex $v_i$ and call it **didegree** of $v_i$, representing respectively in-degree and out-degree of $v_i$, where, $m_i = \delta^-(v_i)$ and $n_i = \delta^+(v_i)$.

**Definition 1.1:** A **diregular tournament** is a diregular complete digraph. (It is clear to see that for the tournament to be diregular the cardinality of its vertex set must be odd.). Thus, a diregular tournament on $(2n+1)$ points is a digraph with in-out-degree-sequence equal to $((n,n),(n,n),\cdots,(n,n))$.

**Definition 1.1:** A **leading diregular tournament**, $G$, on $(2n+1)$ points is the one with the adjacency matrix, $A(G) = [a_{ij}]$, defined as follows:
1) $a_{ij} = 1$ for $(i,j) = (i, i+1), (i, i+2) \cdots, (i, i+n)$ and $i = 1, 2, \cdots, n+1$
2) $a_{ij} = 1$ for
$(i,j) = (n+2,1), (n+3,1), (n+3,2), \cdots, (2n+1,1), (2n+1,2), \cdots, (2n+1,n)$.

**Example:** We state below the adjacency matrix of the **leading** diregular tournament on 9 points:

$$\begin{pmatrix} 0 & 1 & 1 & 1 & 1 & 0 & 0 & 0 & 0 \\ 0 & 0 & 1 & 1 & 1 & 1 & 0 & 0 & 0 \\ 0 & 0 & 0 & 1 & 1 & 1 & 1 & 0 & 0 \\ 0 & 0 & 0 & 0 & 1 & 1 & 1 & 1 & 0 \\ 0 & 0 & 0 & 0 & 0 & 1 & 1 & 1 & 1 \\ 1 & 0 & 0 & 0 & 0 & 0 & 1 & 1 & 1 \\ 1 & 1 & 0 & 0 & 0 & 0 & 0 & 1 & 1 \\ 1 & 1 & 1 & 0 & 0 & 0 & 0 & 0 & 1 \\ 1 & 1 & 1 & 1 & 0 & 0 & 0 & 0 & 0 \end{pmatrix}$$

To explain our procedure of the construction of directed Hamilton circuits, which when taken together provide the desired decomposition, we discuss it first for the special case, i.e., for a leading diregular tournament on $2n+1 = 7$ points. We construct $n = 3$ labeled Hamilton circuits that tightly pack the leading diregular tournament on $2n+1 = 7$ points, i.e., each made up from distinct labeled directed edges such that all these directed labeled edges taken together constitute a labeled leading diregular tournament on $2n+1= 7$ points. For this to achieve we use the existence of $n =3$ coprime numbers $\leq 3$. Since $2n+1 = 7$ is a prime number, such numbers exist. We have $\{1, 2, 3\}$ as the desired set of numbers. We form sequences of the numbers using the numbers in this set as the increment for the sequences, and further by taking the numbers in these sequences as labels of vertices which are joined by the directed edges in the same order, we form the desired Hamilton circuits as follows:

Here, in these sequences the numbers (which are to be taken as suffixes of the vertex labels) are differing from next number (next vertex label) by 1mod (7) for the first sequence, by 2mod (7) for the second sequence, and by 3mod (7) for the third and the last sequence. Thus, we construct three sequences and their associated Hamilton circuits as follows:

**Sequence 1:** (1, 2, 3, 4, 5, 6, 7, 1)
**Hamilton Circuit 1:** $v_1 \to v_2 \to v_3 \to v_4 \to v_5 \to v_6 \to v_7 \to v_1$
**Sequence 2:** (1, 3, 5, 7, 2, 4, 6, 1)

**Hamilton Circuit 2:** $v_1 \to v_3 \to v_5 \to v_7 \to v_2 \to v_4 \to v_6 \to v_1$
**Sequence 3:** (1, 4, 7, 3, 6, 2, 5, 1)
**Hamilton Circuit 3:** $v_1 \to v_4 \to v_7 \to v_3 \to v_6 \to v_2 \to v_5 \to v_1$

**Remark 1.1:** From the above procedure it is the availability of $n = 3$ numbers $\leq n$, namely, {1, 2, 3}, **co-prime to number $(2n+1) = 7$,** which is **important!**

**Remark 1.2:** It is important to see that when the above circuits are taken together they all contain distinct directed edges and together constitute the **leading** diregular tournament on $2n+1 = 7$ points. Thus, these circuits tightly pack the **leading** diregular tournament on $2n+1 = 7$ points by Hamilton circuits as desired! Check the equal indegree and outdegree for every vertex!!

**2. Hamilton Decompositions:** Our idea to prove the existence of decomposition, in terms of Hamilton circuits, for an unlabeled leading diregular tournament is to construct a tight packing for a labeled copy by certain labeled Hamilton circuits which completely define the adjacencies for all the outgoing and incoming directed edges at each labeled vertex.

**Algorithm for construction of labeled Hamilton circuits which tightly pack into a labeled leading diregular tournament on (2n+1) points when (2n+1) is a prime number:**
1) Take $n$ co-prime numbers to given prime number $(2n+1)$. They are $\alpha_1 = 1, \alpha_2 = 2, \cdots, \alpha_n = n$.
2) Construct $n$ number of sequences of numbers such that each sequence starts with number 1 and ends with number 1 and contains in between a particular arrangement of $2n$ numbers from set of numbers {2, 3, ..., 2n+1} as follows:
   **Sequence 1:** (1, $1+\alpha_1 \mod(2n+1)$, $1+2\alpha_1 \mod(2n+1)$, ..., 1)
   **Sequence 2:** (1, $1+\alpha_2 \mod(2n+1)$, $1+2\alpha_2 \mod(2n+1)$, ..., 1)
   .
   .
   .
   **Sequence n:** (1, $1+\alpha_n \mod(2n+1)$, $1+2\alpha_n \mod(2n+1)$, ..., 1)
3) Using the numbers in each sequence as suffixes of the vertices construct the following Hamilton circuits:
   **Hamilton Circuit 1:**

$$v_1 \to v_{(1+\alpha_1 \bmod(2n+1))} \to v_{(1+2\alpha_1 \bmod(2n+1))} \to \cdots \to v_1$$

**Hamilton Circuit 2:**
$$v_1 \to v_{(1+\alpha_2 \bmod(2n+1))} \to v_{(1+2\alpha_2 \bmod(2n+1))} \to \cdots \to v_1$$

$$\vdots$$

**Hamilton Circuit n:**
$$v_1 \to v_{(1+\alpha_n \bmod(2n+1))} \to v_{(1+2\alpha_n \bmod(2n+1))} \to \cdots \to v_1$$

□

**Remark 2.1:** Since $(2n+1)$ is prime, therefore, the numbers $\{1, 2, \ldots, n\}$ are co-prime to $(2n+1)$ and we can build the Hamilton circuits by taking these numbers as $\alpha_1, \alpha_2, \cdots, \alpha_n$ which will tightly pack into leading diregular tournament on $(2n+1)$ points.

We proceed to see that for the **leading diregular tournament on 9 points** this method fails because 9 is not a prime number and so there do not exist 4 co-prime numbers ≤ 4. Using the adjacency matrix given above for the leading tournament on 9 points and using our algorithm to generate Hamiltonian circuits for three numbers {1, 2, 4} co-prime to 9 we have

**Sequence 1:** (1, 2, 3, 4, 5, 6, 7, 8, 9, 1)
**Hamilton Circuit 1:**
$v_1 \to v_2 \to v_3 \to v_4 \to v_5 \to v_6 \to v_7 \to v_8 \to v_9 \to v_1$
**Sequence 2:** (1, 3, 5, 7, 9, 2, 4, 6, 8, 1)
**Hamilton Circuit 2:**
$v_1 \to v_3 \to v_5 \to v_7 \to v_9 \to v_2 \to v_4 \to v_6 \to v_8 \to v_1$
**Sequence 3:** (1, 5, 9, 4, 8, 3, 7, 2, 6, 1)
**Hamilton Circuit 3:**
$v_1 \to v_5 \to v_9 \to v_4 \to v_8 \to v_3 \to v_7 \to v_2 \to v_6 \to v_1$
Now, if we remove these three Hamilton circuits from the leading diregular tournament on 9 points we can check that we will be left with three disjoint triangles (disjoint triangular directed circuits), namely,
$v_1 \to v_4 \to v_7 \to v_1$
$v_2 \to v_5 \to v_8 \to v_2$
$v_3 \to v_6 \to v_9 \to v_3$

Thus, a leading diregular tournament on 9 points cannot be decomposed in Hamiltonian circuits by this method.

We give below two figures. FIG.1 presents the leading regular tournament on 9 points, and FIG.2 depicts its decomposition in terms of 3 disjoint directed Hamilton circuits and 3 disjoint directed and mutually disconnected triangles.

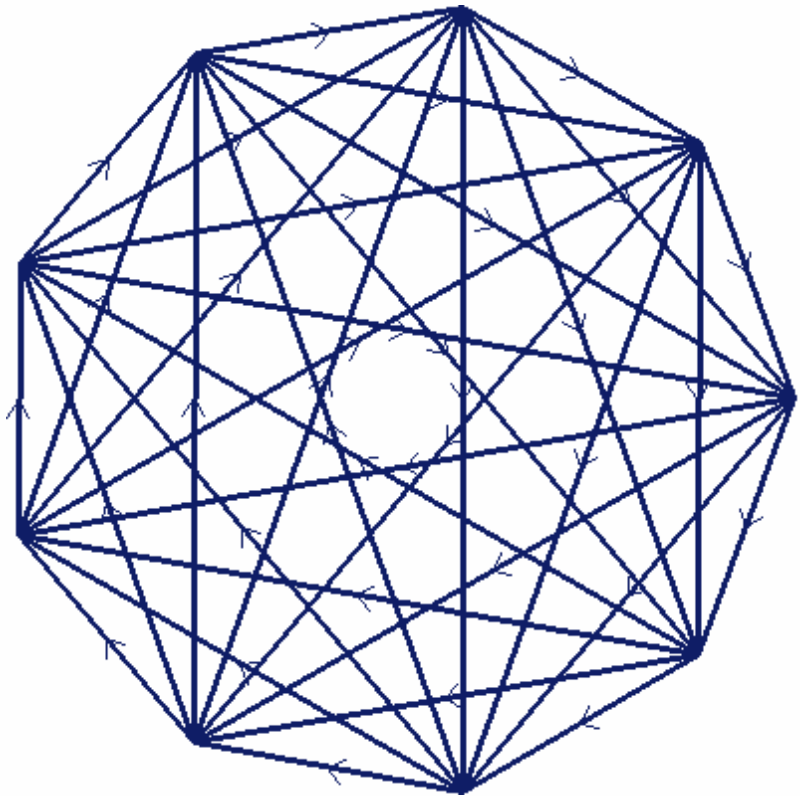

FIG.1

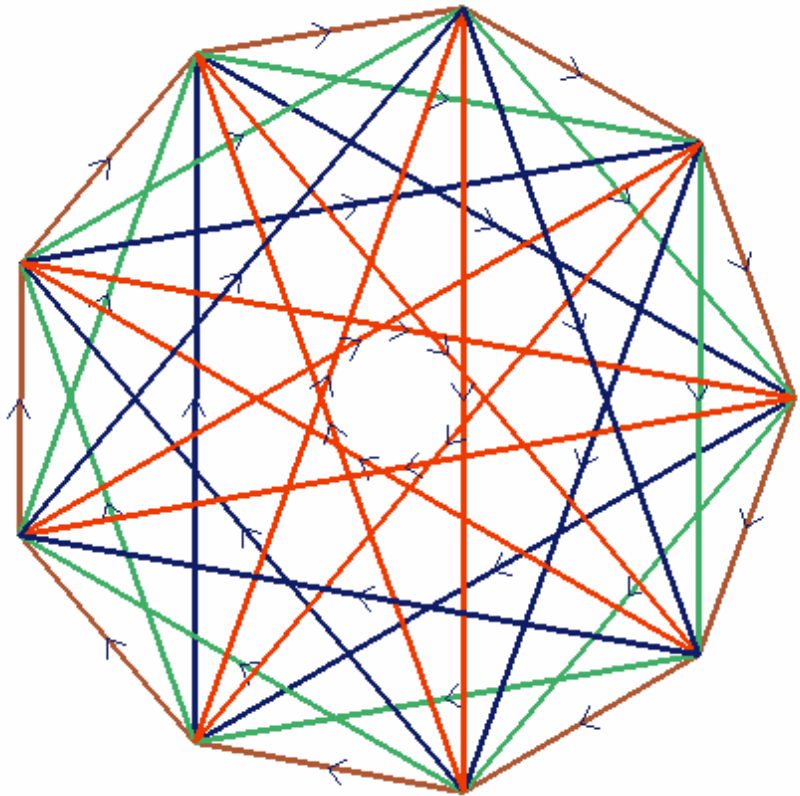

FIG.2

FIG.1 with blue color to all directed edges represents the so called leading tournament on 9 points.
In this FIG.2 we can clearly see the decomposition of this tournament, not in 4 disjoint directed Hamilton circuits but in 3 disjoint directed Hamilton circuits made up of directed edges with brown, green, and orange color respectively, and 3 disjoint directed and mutually disconnected triangles in blue color.

The unlabeled digraph in FIG.1, FIG.2 can be viewed as a **"geometrical objects"** as described below:
We think of FIG.1 as made up of 9 points (vertices) placed equidistant along a circle.
  i) The directed edges joining the neighboring points (shown in brown color in FIG.2) will be said to be at a (hypothetical) distance of exactly "1 unit" from the circle. We call such directed edge the edge of "type 1".

ii) The directed edges joining the next neighboring points (shown in green color in FIG.2) will be said to be at a (hypothetical) distance of exactly "2 units" from the circle. We call such directed edge the edge of "type 2".

iii) The directed edges joining the next neighboring points (shown in blue color in FIG.2) will be said to be at a (hypothetical) distance of exactly "3 units" from the circle. We call such directed edge the edge of "type 3".

iv) The directed edges joining the next neighboring points (shown in orange color in FIG.2) will be said to be at a (hypothetical) distance of exactly "4 units" from the circle. We call such directed edge the edge of "type 4".

Thus, from this classification every directed edge belongs to some one type. Also, there are exactly $2n$ edges of each type. Also, there are exactly two directed edges, one incoming and one outgoing, of each type from every vertex.

It is easy to see that the directed edges at a (hypothetical) distance of exactly "1 unit", "2 units", and "4 units" form the three directed circuits of maximal length (Hamilton circuits) and the directed edges at a (hypothetical) distance of exactly "4 units" form three directed and mutually disconnected triangular circuits.

So, the decomposition for this tournament will be necessarily made up of Hamilton circuits with edges belonging to mixed types, i.e. belonging to different (hypothetical) distance types.

We now proceed with the formal proof of our result:

**Theorem 2.1:** Every leading diregular tournament on $(2n+1)$ points can be decomposed in $n$ directed Hamilton circuits if $(2n+1)$ is a prime number.

**Proof:** When $(2n+1)$ is prime then the first $n$ numbers $\{1, 2, \ldots, n\}$ are always co-prime to $(2n+1)$. So, we can always construct $n$ distinct sequences as described in the algorithm. From those $n$ distinct sequences, starting and ending with the number 1, we can construct $n$ distinct (disjoint) Hamilton circuits, each containing distinct directed edges, and all these edges taken together clearly form the leading diregular tournament, i.e. these distinct Hamilton circuits tightly pack in the leading diregular tournament, since for every vertex there is exactly one incoming directed edge and exactly one outgoing directed edge, each such vertex appear in every directed circuit and these incoming and outgoing edges in each Hamilton circuit are distinct. So, every vertex has exactly $n$ incoming edges and exactly $n$ outgoing edges. Hence etc.

□

There can exist certain tournaments on $(2n+1)$ points, belonging **some isomorphism type**, for which the Hamilton decomposition can exist for all $n$ belonging to set $\{1, 2, \ldots\}$. For the demonstration of the existence of Hamilton decomposition for these tournaments we discuss a geometrical construction as is provided for complete graphs in [2] for a type of tournament that automatically builds up in the process of its construction defined below.

**Theorem 2.2:** There exists a diregular tournament on $(2n+1)$ points of certain isomorphism type, i.e. belonging to some certain isomorphism class, which can be decomposed into $n$ disjoint directed Hamilton circuits where $(2n+1)$ is not necessarily a prime number.

**Proof:** A diregular tournament on $(2n+1)$ vertices has $n(2n+1)$ directed edges. A directed Hamilton circuit on $(2n+1)$ points contains $(2n+1)$ directed edges. Therefore, the number of edge disjoint directed Hamilton circuits cannot exceed $n$.
To show that there are exactly $n$ such directed circuits which tightly pack certain diregular tournament on $(2n+1)$ points we proceed with the following geometrical construction:
   1) Draw a circle (with a dotted line showing its boundary) with center taken as a vertex with label $v_1$.
   2) Draw a horizontally going diameter (along x-axis) cutting the dotted circle at two points and take the points where this diameter cuts the dotted circle as vertices with labels $v_2$ (on the negative side of the x-axis) and $v_{(2n+1)}$ (on the positive side of the x-axis) respectively. Give direction to these edges as $v_1 \to v_2$ (i.e. $v_1$ to $v_2$) and $v_{(2n+1)} \to v_1$ ( i.e. $v_{(2n+1)}$ to $v_1$) respectively.

3) Draw diameters, each one cutting the dotted circle at two points, making angles $\left(\frac{\pi}{n}\right), \left(\frac{2\pi}{n}\right), \cdots, \left(\frac{(n-1)\pi}{n}\right)$ with the horizontal measured in the clockwise direction, and take the points where these diameters cut the dotted circle as vertices with serially chosen labels $v_3, v_5, \cdots, v_{(2n-1)}$ to points taken as vertices on the circle by starting at vertex with label $v_2$ and moving along the circle in clockwise direction, and as vertices with labels $v_4, v_6, \cdots, v_{2n}$ by starting at vertex with label $v_2$ and moving along the circle in anticlockwise direction.

4) Create directed edges: $v_2 \to v_3$ (i.e. $v_2$ to $v_3$), $v_3 \to v_4$, $v_4 \to v_5$, ...., $v_{2n} \to v_{(2n+1)}$. These directed edges along with the directed edges created in the second step will constitute the directed Hamilton circuit: $v_1 \to v_2 \to v_3 \to \cdots \to v_{2n} \to v_{(2n+1)} \to v_1$.

5) Now Keep the vertices along with their labels fixed on the circle and rotate the geometrical pattern formed by the directed edges in a sequence respectively by angles $\left(\frac{\pi}{n}\right), \left(\frac{2\pi}{n}\right), \cdots, \left(\frac{(n-1)\pi}{n}\right)$ in succession, each time recording the new disjoint Hamilton circuit that gets created at each rotation, like, $v_1 \to v_3 \to v_5 \to v_2 \to \cdots \to v_{2n} \to v_1$, that is formed after first rotation, namely, by $\left(\frac{\pi}{n}\right)$, etc.

Note that each rotation produces a directed Hamilton circuit that has no directed edge common with any of the previous ones. Thus we have in all $n$ disjoint directed Hamilton circuits. Now, also for each vertex we have one incoming directed edge and one outgoing directed edge and these edges are different in different circuits. Thus, for every vertex we have exactly $n$ incoming directed edges and exactly $n$ outgoing directed edges. Thus, these Hamilton circuits form (or tightly pack into) certain diregular tournament whose adjacency relations get defined in its construction procedure, as desired! Hence etc. □

**Conclusion:** Leading diregular tournament can be decomposed into directed Hamilton circuits with edges in any Hamilton circuit belong to a fixed distance type when (2n+1) is prime. Some type of tournaments can always be decomposed in Hamilton circuits.

## Acknowledgements

Many thanks are due to Dr. M. R. Modak and Mr. Riko Winterle for encouragements. Many thanks are also due to Mr. Abhijit Patwardhan for providing the colored figures for the leading diregular tournament on 9 points.